\documentclass[12pt,times new roman]{article}
\usepackage{amsfonts}
\usepackage{amssymb}
\usepackage{amsmath}
\usepackage{amscd}
\usepackage{times}
\usepackage{amsthm}
\usepackage{hyperref}
\usepackage[margin=1.5in]{geometry}
\usepackage{multirow}
\usepackage{indentfirst}
\usepackage{dsfont}
\usepackage{authblk}
\usepackage{setspace}
\newtheorem{Theorem}{Theorem}[section]
\newtheorem{Lemma}[Theorem]{Lemma}
\newtheorem{proposition}[Theorem]{Proposition}
\newtheorem{innercustomthm}{Theorem}
\newtheorem{innercustomprop}{Proposition}

\newenvironment{mythm}[1]
  {\renewcommand\theinnercustomthm{#1}\innercustomthm}
  {\endinnercustomthm}
\newenvironment{myprop}[1]
  {\renewcommand\theinnercustomprop{#1}\innercustomprop}
  {\endinnercustomprop}
\title{$G_0$ of affine, simplicial toric varieties}

\author{Zeyu Shen}
\date{June 3rd, 2024}

\begin{document}
\maketitle

\begin{abstract}
Let $X$ be an affine, simplicial toric variety 
over a field. 
Let $G_0$ denote the Grothendieck group of coherent sheaves 
on a Noetherian scheme and let $F^1G_0$ denote the first 
step of the filtration on $G_0$ by codimension of support.
Then $G_0(X)\cong\mathbb{Z}\oplus F^1G_0(X)$ and 
$F^1G_0(X)$ is a finite abelian group.
In dimension 2, we show that $F^1G_0(X)$ is a finite cyclic 
group and determine its order.
In dimension 3, $F^1G_0(X)$ is determined up to a group 
extension of the Chow group $A^1(X)$ by the Chow group 
$A^2(X)$.
We determine the order of the Chow group $A^1(X)$ in this 
case.
A conjecture on the orders of $A^1(X)$ and $A^2(X)$ is
formulated for all dimensions.\\
\end{abstract}
\section{Introduction}
Several authors such as Morelli, Gubeladze, Cortiñas et 
al. have studied the algebraic $K$-theory of toric varieties 
in \cite{M},\cite{G2},\cite{CHWW}.
Others, such as Joshua and Krishna, have considered the 
equivariant algebraic $G$-theory of toric stacks in 
\cite{JK}.
However, there have been no results specific to the 
algebraic $G$-theory of affine, simplicial toric varieties.
We seek to address this issue by determining the structure 
of $G_0$, the Grothendieck group of coherent sheaves.\par
Let $X$ be an affine, simplicial toric variety over a
field $k$.
If $X$ has dimension zero, then $X\cong\operatorname{Spec}
(k)$, so that $G_0(X)\cong\mathbb{Z}$.
If $X$ has dimension one, then $X\cong\mathbb{A}_k^1$, so 
that $G_0(X)\cong\mathbb{Z}$ as well.
We show that $G_0(X)\cong\mathbb{Z}\oplus\mathbb{Z}/\delta\mathbb{Z}$ 
when $X$ has dimension 2, where $\delta$ is the determinant 
of the matrix taking the minimal generators of the fan of 
$X$ as its columns. 
In dimension 3, the Brown-Gersten-Quillen spectral sequence 
and Fulton's Riemann-Roch theorem for algebraic schemes 
together determine $F^1G_0(X)$ as a group extension of the 
Chow group $A^1(X)$ by the Chow group $A^2(X)$.\
Here $F^1G_0(X)$ is the first step of the filtration on 
$G_0(X)$ by codimension of support. 
By simplifying a general 3-dimensional simplicial cone
in $\mathbb{R}^3$ to a nice form, we compute the order of 
the Chow group $A^1(X)$.\par
We first prove the following general result:
\begin{myprop}{2.1}
Let $X$ be an affine, simplicial toric variety over a
field $k$.\\
Then $G_0(X)\cong\mathbb{Z}\oplus F^1G_0(X)$ and $F^1G_0(X)$ 
is a finite abelian group.\par
\end{myprop}
\noindent
Here is our result in dimension 2:
\begin{mythm}{3.3}
Let $X$ be an affine toric surface over a field $k$ of the 
form 
$\operatorname{Spec}(k[\sigma^{\vee}\cap\mathbb{Z}^2])$, 
where $\sigma$ is a 2-dimensional strongly convex, rational, 
polyhedral cone in $\mathbb{R}^2$.
Then 
$G_0(X)\cong\mathbb{Z}\oplus\mathbb{Z}/\delta\mathbb{Z}$, 
where $\delta$ is the determinant of the matrix taking the 
minimal generators of the cone $\sigma$ as its columns.\\    
\end{mythm}
\noindent
In dimension 3, our result is Theorem 4.2:\\
\begin{mythm}{4.2}
Let $X$ be an affine, simplicial toric 3-fold over a field 
$k$.
Then $F^1G_0(X)$ is an extension of the Chow group $A^1(X)$ 
by the Chow group $A^2(X)$.
\end{mythm}
\noindent
In Theorem 4.4, we compute the order of the Chow group 
$A^1(X)$ of any 3-dimensional affine, simplicial toric 
variety over $k$.\\
\begin{mythm}{4.4}
Let $X=\operatorname{Spec}(k[\sigma^{\vee}\cap\mathbb{Z}^3])$ be an affine, simplicial 
toric 3-fold over a field 
$k$ associated to a 3-dimensional simplicial cone $\sigma$ 
in $\mathbb{R}^3$.
Then the Chow group $A^1(X)$ has order $|\delta|$, where 
$\delta$ is the determinant of the matrix taking the minimal 
generators of the cone $\sigma$ as its columns.
\end{mythm}
The paper is organized as follows.
In the preliminary section, we prove the general result 
Proposition 2.1.
In section 3, we prove that any 2-dimensional strongly 
convex, rational, polyhedral cone in $\mathbb{R}^2$ can be 
mapped via a transformation in $GL(2,\mathbb{Z})$ to a cone 
in $\mathbb{R}^2$ with minimal generators of the form 
$e_1,ae_1+be_2$, $a,b\in\mathbb{Z}$.
We also observe that any transformation in 
$GL(n,\mathbb{Z})$ 
which maps an $n$-dimensional simplicial cone $\sigma$ to 
another such cone $\tau$ induces a ring isomorphism\\
$k[\sigma^{\vee}\cap\mathbb{Z}^n]\xrightarrow{\cong}
k[\tau^{\vee}\cap\mathbb{Z}^n]$.\par
In section 3 we prove Theorem 3.3.
In section 4, we prove Theorems 4.2 and 4.4.
In section 5, we formulate a conjecture on the orders of the 
Chow groups $A^1(X)$ and $A^2(X)$.\par
{\bf Acknowledgements.} 
The author would like to thank his 
advisor Professor Charles Weibel for continual guidance and 
helpful comments throughout the preparation of this paper.
The author is also grateful for being supported by NSF Grant 
DMS 2001417.\par
\section{Preliminaries}
\noindent
Let $X$ be an affine, simplicial toric variety over a field 
$k$.
Let $F^iG_0(X)$ be the $i$-th step of the filtration on 
$G_0(X)$ by codimension of support.
Let $E$ denote the Brown-Gersten-Quillen spectral sequence 
for the Noetherian scheme $X$.
\begin{proposition}
Let $X$ be an affine, simplicial toric variety over a 
field $k$.
Then $G_0(X)\cong\mathbb{Z}\oplus F^1G_0(X)$,
and $F^1G_0(X)$ is a finite abelian group.
And $F^nG_0(X)=0$, where $n$ is the Krull dimension of $X$.
\end{proposition}
\begin{proof}
By construction, we have $E_{\infty}^{0,0}\cong 
E_2^{0,0}=A^0(X)\cong\mathbb{Z}$.
And we have an exact sequence
\begin{equation*}
0\rightarrow F^1G_0(X)\rightarrow G_0(X)\rightarrow 
G_0(X)/F^1G_0(X)\rightarrow 0.
\end{equation*}
By construction, we have $G_0(X)/F^1G_0(X)\cong E_{\infty}^{0,0}$, so $G_0(X)/F^1G_0(X)\cong\mathbb{Z}$.\\
Thus the exact sequence above splits to yield
$G_0(X)\cong\mathbb{Z}\oplus\ F^1G_0(X)$.\par
Let $n$ be the Krull dimension of $X$.
If $n=0$, then $F^1G_0(X)=0$ by construction.
We may assume $n>0$.
By construction, $F^{n+1}G_0(X)=0$.
Since $X$ is a non-complete toric variety, $A^n(X)=A_0(X)=0$ 
follows from the discussion after Example 
2.3 in \cite{FS}. 
So $E_{\infty}^{n,-n}\cong E_2^{n,-n}=A^n(X)=0$.
Thus, 
\begin{equation*}
F^nG_0(X)=F^nG_0(X)/F^{n+1}G_0(X)\cong E_{\infty}^{n,-n}\cong 
A^n(X)=0.
\end{equation*}
Consider the exact sequences 
\begin{equation}\label{1}
    \begin{aligned}
0\rightarrow F^1G_0(X)\rightarrow G_0(X)\rightarrow 
G_0(X)/F^1G_0(X)\rightarrow 0    \\
0\rightarrow F^2G_0(X)\rightarrow F^1G_0(X)\rightarrow 
F^1G_0(X)/F^2G_0(X)\rightarrow 0 \\
...\\
0\rightarrow F^{n-1}G_0(X)\rightarrow F^{n-
2}G_0(X)\rightarrow F^{n-2}G_0(X)/F^{n-1}G_0(X)\rightarrow 0
  \end{aligned}
\end{equation}
By Fulton's Riemann-Roch Theorem for Algebraic Schemes 
\cite{F} 
we have \\
$G_0(X)\otimes\mathbb{Q}\cong A^*(X)\otimes\mathbb{Q}$.
By a result of Gubeladze in \cite{G1},
$K_0(X)\cong\mathbb{Z}$ since $X$ is an affine toric 
variety.
Since $X$ is a simplicial toric variety, the corollary of 
section 2.2 in \cite{BV} states that
the canonical homomorphism 
$K_0(X)\otimes\mathbb{Q}\rightarrow G_0(X)\otimes\mathbb{Q}$
is surjective.
Hence, $G_0(X)\otimes\mathbb{Q}$ has dimension at most 1.
Since $G_0(X)$ contains $\mathbb{Z}$ as a direct summand, 
$G_0(X)\otimes\mathbb{Q}$ has dimension at least 1.
Therefore, $G_0(X)\otimes\mathbb{Q}$ has dimension 1.
So $A^*(X)\otimes\mathbb{Q}$ also 
has dimension 1.
Since $A^0(X)\otimes\mathbb{Q}=\mathbb{Q}$, the sum of 
dimensions of $A^i(X)\otimes\mathbb{Q}$ for all $i>0$ is 
zero.
By the lemma in \cite[12.5.1]{CLS} all the Chow groups 
$A^p(X)$ of the toric variety $X$ are finitely generated 
abelian groups.
So the Chow groups $A^i(X)$ are 
all finite for $i>0$.
Since $F^iG_0(X)/F^{i+1}G_0(X)\cong E_{\infty}^{i,-
i}$ for every $i$, and that $E_{\infty}^{i,-i}$ is always a 
subquotient of the abelian group $E_2^{i,-i}=A^i(X)$ for 
every $i$, we deduce from the exact sequences \eqref{1}
that $F^1G_0(X)$ is a finite abelian group.
\end{proof}
\begin{Lemma}
Let $\sigma$ be an $n$-dimensional strongly convex, rational, 
polyhedral cone in $\mathbb{R}^n$. 
Let $R=k[\sigma^{\vee}\cap\mathbb{Z}^n]$.
Then $R^{\times}=k^{\times}$.
\end{Lemma}
\begin{proof}
Since $k$ is a subring of $R$, $k^{\times}\subseteq R^{\times}$.
Since the polyhedral cone $\sigma^{\vee}$ is strongly 
convex, it does not contain any positive dimensional vector 
subspace of $\mathbb{R}^n$.
Since multiplication in the ring $R$ corresponds to vector 
addition in $\mathbb{R}^n$, we conclude that 
$R^{\times}\subseteq k^{\times}$.
Otherwise, there is a nonzero vector 
$v\in\sigma^{\vee}\cap\mathbb{Z}^n$ such that $-
v\in\sigma^{\vee}\cap\mathbb{Z}^n$.
But this contradicts the hypothesis that $\sigma^{\vee}$ is 
does not contain any positive dimensional vector subspace of 
$\mathbb{R}^n$.
Hence, we have $R^{\times}=k^{\times}$.
\end{proof}
\section{dim($X$)=2}
\begin{Lemma}
Let $\sigma$ be a 2-dimensional strongly convex, rational, 
polyhedral cone in $\mathbb{R}^2$. There exists a matrix 
$M\in GL(2,\mathbb{Z})$ such that left multiplication by $M$ 
maps the cone $\sigma$ isomorphically onto
the cone in $\mathbb{R}^2$ generated by $e_1,ae_1+be_2$, 
where $a,b$ are relatively prime integers and $a,b>0$.\\
\end{Lemma}
\begin{proof}
Let $(x_1,y_1)$, $(x_2,y_2)\in \mathbb{Z}^2$ be the minimal 
generators for the cone $\sigma$.
Since the integers $x_1,y_1$ are relatively prime, there 
exist integers 
$a,b$ such that $ax_1+by_1=1$.\par
Let $A$ be the matrix 
$\begin{pmatrix}
a & b\\
-y_1 & x_1\\
\end{pmatrix}.$
Let $B$ be the matrix 
$\begin{pmatrix}
1 & 0\\
0 & -1
\end{pmatrix}.$
Since the matrix $A$ has determinant $ax_1-b(-y_1)=1$,
$A$ and $B$ are in $GL(2,\mathbb{Z})$.
Note that left multiplication by $A$ maps $(x_1,y_1)$ to 
$e_1$ and maps the minimal generator 
$(x_2,y_2)$ of the cone $\sigma$ to $(ax_2+by_2,-x_2y_1+x_1y_2)$.
Note that the second coordinate  $t=x_1y_2-
x_2y_1$ must be nonzero, since the vectors $(x_1,y_1),
(x_2,y_2)$ are minimal generators for the 2-dimensional cone 
$\sigma$, which are linearly independent over 
$\mathbb{R}$.\par
If $t$ is positive, left multiplication by $A\in 
GL(2,\mathbb{Z})$ maps the cone $\sigma$ isomorphically onto 
the cone in 
$\mathbb{R}^2$ generated by $e_1,pe_1+qe_2$ for some 
integers $p,q$ and $q>0$.
If $t$ is negative, then left multiplication by 
$BA\in GL(2,\mathbb{Z})$ maps the cone $\sigma$ isomorphically onto 
the cone in $\mathbb{R}^2$ generated by $e_1,pe_1+qe_2$ for some integers 
$p,q$ and $q>0$.
Since the ray generated by $pe_1+qe_2$ is the same as the ray generated by 
$p'e_1+q'e_2$, where $p'=\frac{p}{d}$, $q'=\frac{q}{d}$ and $d=gcd(p,q)$,
we may assume that the integers $p,q$ are relatively 
prime.\par
Let $m$ be an integer and let $C$ be the matrix 
$\begin{pmatrix}
1 & m\\
0 & 1
\end{pmatrix}$ 
in $GL(2,\mathbb{Z})$.
Note that left multiplication by $C$ maps the vector $(p,q)$ 
to $(p+qm,q)$.
Choose $m$ such that $p+qm>0$. 
If $x_1y_2-x_2y_1>0$, let $M=CA$.
If $x_1y_2-x_2y_1<0$, let $M=CBA$.
Then in both cases $M\in GL(2,\mathbb{Z})$ and left 
multiplication by $M$ maps the cone $\sigma$ isomorphically 
onto the cone in $\mathbb{R}^2$ generated by 
$e_1,ae_1+be_2$, where $a,b>0$ are relatively prime integers.
\end{proof}
\begin{Lemma}
Let $A\in GL(n,\mathbb{Z})$ be the matrix such that 
$Au_i=v_i$, where $u_1,u_2,...,u_n$ are the minimal 
generators of an $n$-dimensional simplicial cone $\sigma$ in 
$\mathbb{R}^n$ and $v_1,v_2,...,v_n$ are the minimal 
generators of an $n$-dimensional simplicial cone $\tau$ in 
$\mathbb{R}^n$.\\
Then left multiplication by the matrix $(A^{-1})^{t}$ 
induces a ring isomorphism\\
$k[\sigma^{\vee}\cap\mathbb{Z}^n]\xrightarrow{\cong} 
k[\tau^{\vee}\cap\mathbb{Z}^n]$.\\
\end{Lemma}
\begin{proof}
Let $x\in\sigma^{\vee}$ and $y\in\tau$.
By construction of the matrix $A$, left multiplication by 
$A^{-1}$ maps $y$ into the cone $\sigma$.
Hence, $\langle (A^{-1})^{t}x,y\rangle =x^{t}A^{-1}y=\langle x,A^{-1}y\rangle \geq 0$,
i.e., left multiplication by $(A^{-1})^{t}$ maps the cone 
$\sigma^{\vee}$ into the cone $\tau^{\vee}$.
Let $w\in\tau^{\vee}$.
Set $v=A^tw$.
Then $v\in\mathbb{R}^n$ satisfies $\langle v,u\rangle=\langle A^tw,u\rangle =w^{t}Au=
\langle w,Au\rangle \geq 0$ for every $u\in\sigma$, since left 
multiplication by $A$ maps $\sigma$ into $\tau$,
i.e., left multiplication by $(A^{-1})^{t}$ is a bijection 
from the cone $\sigma^{\vee}$ to the cone $\tau^{\vee}$.
Since the matrix $A$ is in $GL(n,\mathbb{Z})$, so is the 
matrix $(A^{-1})^{t}$.
Hence, left multiplication by $(A^{-1})^{t}$ maps the 
semigroup $\sigma^{\vee}\cap\mathbb{Z}^n$ into the semigroup
$\tau^{\vee}\cap\mathbb{Z}^n$.
And this restriction of left multiplication by $(A^{-1})^{t}$
is a semigroup 
isomorphism:
$\sigma^{\vee}\cap\mathbb{Z}^n\xrightarrow{\cong}
\tau^{\vee}\cap\mathbb{Z}^n$.
Applying the semigroup algebra functor over the field $k$, 
we obtain an induced ring isomorphism:
$k[\sigma^{\vee}\cap\mathbb{Z}^n]\xrightarrow{\cong}k[\tau^{\vee}\cap\mathbb{Z}^n]$.
\end{proof}
\begin{Theorem}
Let $X=\operatorname{Spec}(k[\sigma^{\vee}\cap\mathbb{Z}^2])$
be an affine, simplicial toric surface over a field $k$ 
associated to a 2-dimensional simplicial cone $\sigma$ in 
$\mathbb{R}^2$.
Then 
$G_0(X)\cong\mathbb{Z}\oplus\mathbb{Z}/\delta\mathbb{Z}$, 
where $\delta$ is the determinant of the matrix taking the 
minimal generators of the cone $\sigma$ as its columns.\\
\end{Theorem}
\begin{proof}
Let $E$ denote the Brown-Gersten-Quillen spectral sequence 
for the Noetherian scheme $X$ of dimension 2.
By construction, the differentials entering and leaving the 
diagonal terms $E_2^{0,0}$ and $E_2^{1,-1}$ are both zero.
Hence, $E_{\infty}^{0,0}\cong 
E_2^{0,0}=A^0(X)\cong\mathbb{Z}$ and 
$E_{\infty}^{1,-1}\cong E_2^{1,-1}=A^1(X)$.\par
Since $X$ is a variety, the Chow group $A^1(X)$ is 
isomorphic to the divisor class group $\mathrm{Cl}(X)$.
Now we compute the Chow group $A^1(X)$.
Let $M$ be the character lattice for the toric variety 
$X$.
There is an exact sequence\par
\begin{equation}\label{2}
M\xrightarrow {\alpha} \mathrm{Div}_T(X)\rightarrow \mathrm{Cl}(X)\rightarrow 
0
\end{equation}
where $T$ denotes the torus of the toric variety $X$.
Here $\alpha$ maps every $m\in M$ to the $T$-invariant 
divisor $\mathrm{div}(\chi^{m})$, and the map 
$\mathrm{Div}_T(X)\rightarrow 
\mathrm{Cl}(X)$ maps every $T$-invariant Weil divisor on 
$X$ to its divisor class.
We have $M=\mathbb{Z}^2$ and 
$\mathrm{Div}_T(X)\cong\mathbb{Z}^2$, since there are two rays in the 
fan of the toric variety $X$.
The matrix representing the map $\alpha$ under the standard 
bases takes the minimal generators of the cone $\sigma$ as 
its rows.
By Proposition 2.1 and Lemma 3.1 above, we may assume that 
the minimal generators of the cone $\sigma$ are $e_1, 
ae_1+be_2$, $a,b\in\mathbb{Z}$.\par
Hence, the matrix representing the map $\alpha$ under the 
standard bases is
$\begin{pmatrix}
1 & 0\\
a & b\\
\end{pmatrix}$.
By adding $-a$ times of the first row of the matrix above to 
its second row, we obtain\\
$\begin{pmatrix}
1 & 0\\
0 & b\\
\end{pmatrix}$ as its Smith Normal form.
Hence we deduce that $\mathrm{Cl}
(X)\cong\mathbb{Z}/b\mathbb{Z}$.
Note that by construction, we have $\delta=b$.
Therefore, we have $A^1(X)=\mathrm{Cl}(X)\cong\mathbb{Z}/\delta\mathbb{Z}$.\par
By Proposition 2.1 above, $F^2G_0(X)=0$.
Hence, 
\begin{equation*}
F^1G_0(X)\cong F^1G_0(X)/F^2G_0(X)\cong 
E_{\infty}^{1,-1}\cong\mathbb{Z}/\delta\mathbb{Z}    
\end{equation*}
We conclude that 
$G_0(X)\cong\mathbb{Z}\oplus\mathbb{Z}/\delta\mathbb{Z}$.
\end{proof}
\section{dim($X$)=3}
Let $E_{*}^{*,*}=E_*^{*,*}(X)$ be the Brown-Gersten-Quillen 
spectral sequence of the Noetherian scheme $X$.
Let $A^2(X)$ be the Chow group of codimension 2 cycles on 
$X$.
\begin{proposition}
Let $X$ be an affine, simplicial toric 3-fold over a field 
$k$.
Then $A^2(X)\cong E_2^{2,-2}\cong E_{\infty}^{2,-2}$.\\
\end{proposition}
\begin{proof}
Say $X=\operatorname{Spec}(R)$, where
$R=k[\sigma^{\vee}\cap\mathbb{Z}^3]$ for some 3-
dimensional simplicial cone $\sigma$ in $\mathbb{R}^3$.
By Lemma 2.2 above, $R^{\times}=k^{\times}$.
Since the $k$-algebra generators of $R$ are Laurent 
monomials in the variables $x,y,z$, $R$ is a subring of the 
rational function field $k(x,y,z)$. Hence, $R$ is an 
integral domain.\par
By Exercise \cite[II 6.9]{W}, there is an exact sequence
\begin{equation*}
0\rightarrow R^{\times}\rightarrow F^{\times}\xrightarrow{\Delta} D^1(R)\rightarrow CH^1(R)\rightarrow 0
\end{equation*}
Here $F$ is the field of fractions of $R$ and $D^1(R)$ is 
the free abelian group on the height one prime ideals of $R$.
By Proposition \cite[V 9.2]{W}, the map $\Delta$ above is 
the differential of the Brown-Gersten-Quillen spectral 
sequence $d_1:E_1^{0,-1}\rightarrow E_1^{1,-1}$ for the ring 
$R$.
Thus, the kernel of this differential is the image of the 
natural map $R^{\times}\rightarrow F^{\times}$, which is 
$R^{\times}=k^{\times}$.
For the Brown-Gersten-Quillen spectral sequence $E(R)$ of 
the ring $R$, we have $E_1^{-1,-1}(R)=0$, so
$E_2^{0,-1}(R)=\ker(d_1:E_1^{0,-1}(R)\rightarrow 
E_1^{1,-1}(R))$.
Hence, we have $E_2^{0,-1}(R)=\ker(d_1:E_1^{0,-1}
(R)\rightarrow E_1^{1,-1}(R))=k^{\times}$.\par
The flat ring homomorphism $k\rightarrow 
R$ induces a covariant map of Brown-Gersten-Quillen 
spectral sequences.
Let $E(k)$ denote the Brown-Gersten-Quillen spectral 
sequence of the ring $k$.
We have a commutative diagram coming from the map of 
spectral sequences $E(k)\rightarrow E(R)$:\par
$\begin{CD}
k^{\times}=E_2^{0,-1}(k) @>d_2(k)>> E_2^{2,-2}(k)=0\\
@VVV                        @VVV\\
k^{\times}=E_2^{0,-1}(R) @>d_2(R)>> E_2^{2,-2}(R)=A^2(X).
\end{CD}$

Since $k$ has dimension zero, $E_2^{2,-2}(k)=0$.
The map $E_2^{0,-1}(k)\rightarrow E_2^{0,-1}(R)$ is induced 
by the map $E_1^{0,-1}(k)\rightarrow E_1^{0,-1}(R)$.
By construction, this is the natural map of groups of units 
$k^{\times}\rightarrow F^{\times}$, which is induced by the 
ring map $k\rightarrow F$.
So this map is injective.
Hence, the vertical map $E_2^{0,-1}(k)\rightarrow E_2^{0,-1}
(R)$ is an isomorphism.
By the commutative diagram above, the map $d_2(R)$ is zero.
By construction, $E_2^{4,-3}(R)=0$.
So the differentials coming into and leaving from $A^2(X)$ 
are both zero.
Therefore, $A^2(X)\cong E_2^{2,-2}\cong 
E_{\infty}^{2,-2}$.
\end{proof}
\begin{Theorem}
Let $X$ be an affine, simplicial toric 3-fold over a field 
$k$.
Then $G_0(X)\cong\mathbb{Z}\oplus F^1G_0(X)$, and
$F^1G_0(X)$ is an extension of the Chow group $A^1(X)$ by 
the Chow group $A^2(X)$.
\end{Theorem}
\begin{proof}
By Proposition 2.1 above, we have 
$G_0(X)\cong\mathbb{Z}\oplus F^1G_0(X)$ and $F^3G_0(X)=0$.
We have an exact sequence
\begin{equation}\label{3}
0\rightarrow F^2G_0(X)\rightarrow F^1G_0(X)\rightarrow F^1G_0(X)/F^2G_0(X)\rightarrow 0.  
\end{equation}
We have $F^2G_0(X)\cong 
F^2G_0(X)/F^3G_0(X)\cong E_{\infty}^{2,-2}\cong A^2(X)$ by
Proposition 4.1 above.
Since $F^1G_0(X)/F^2G_0(X)\cong E_{\infty}^{1,-1}\cong 
A^1(X)$, \eqref{3} simplifies to
the exact sequence
\begin{equation*}
0\rightarrow A^2(X)\rightarrow F^1G_0(X)\rightarrow 
A^1(X)\rightarrow 0.\qedhere
\end{equation*}
\end{proof}
\begin{Lemma}
Let $\sigma$ be a 3-dimensional simplicial cone in 
$\mathbb{R}^3$. Then there is some matrix $A\in GL(3,\mathbb{Z})$ such that left multiplication by $A$
transforms $\sigma$ into a 3-dimensional simplicial cone 
with minimal generators of the form 
$e_1,ae_1+be_2,ce_1+de_2+ee_3$, where $a,b,c,d,e$ are all 
integers and $e\neq 0$.\\
\end{Lemma}
\begin{proof}
Let $u_1,u_2,u_3$ be the three minimal generators of the 3-
dimensional simplicial cone $\sigma$.
Since $u_1$ is a unimodular column, it can extended to an 
invertible matrix $M$ over the integers.
Hence, left multiplication by the matrix $M^{-1}$ maps $u_1$ 
to $e_1\in\mathbb{R}^3$. 
Thus we may assume that the cone $\sigma$ has minimal 
generators 
\begin{equation*}
u_1=e_1, u_2=a_1e_1+a_2e_2+a_3e_3, u_3=b_1e_1+b_2e_2+b_3e_3
\end{equation*}
Set $d=gcd(a_2,a_3), A_3=-\frac{a_3}{d}, B_3=\frac{a_2}
{d}$. Then $A_3$, $B_3$ are relatively prime integers.
Hence, there exist integers $A_2,B_2$ such that $A_2B_3-
A_3B_2=1$.
Let $\phi:\mathbb{R}^3\rightarrow\mathbb{R}^3$ be the 
$\mathbb{R}$-linear map defined by mapping $e_1$ to $e_1$,
$e_2$ to $A_1e_1+A_2e_2+A_3e_3$ and $e_3$ to 
$B_1e_1+B_2e_2+B_3e_3$, where the $A_i$'s and the $B_i$'s 
are all integers.
Then the coefficient of $e_3$ in 
$\phi(u_2)$ is $a_2A_3+a_3B_3$, which is 
zero, by the construction of the integers $A_3,B_3$.
Since $A_2B_3-A_3B_2=1$, the determinant of the linear map 
$\phi$ is 1.
Hence, $\phi$ is in $SL(3,\mathbb{Z})$ and it transforms 
the 3-dimensional simplicial cone $\tau$ to another 3-
dimensional simplicial cone, whose minimal generators have 
the form $e_1,ae_1+be_2,ce_1+de_2+ee_3$, where $a,b,c,d,e$ 
are integers and $e\neq 0$.\\
Let $A=NM^{-1}$, where $N$ is the matrix representing $\phi$
under the standard basis of $\mathbb{R}^3$.
Then $A\in GL(3,\mathbb{Z})$ is the matrix we are looking 
for.
\end{proof}
\begin{Theorem}
Let $X=\operatorname{Spec}(k[\sigma^{\vee}\cap\mathbb{Z}^3])$ be an affine, simplicial 
toric 3-fold over a field 
$k$ associated to a 3-dimensional simplicial cone $\sigma$ 
in $\mathbb{R}^3$.
Then the Chow group $A^1(X)$ has order $|\delta|$, where 
$\delta$ is the determinant of the matrix taking the minimal 
generators of the cone $\sigma$ as its columns.\\
\end{Theorem}
\begin{proof}
Note that the Chow group $A^1(X)$ is the divisor class group 
$\mathrm{Cl}(X)$.
As noted in the proof of Theorem 3.3 \ref{2} there is an 
exact sequence
\begin{equation*}
M\xrightarrow{\alpha}\mathrm{Div}_T(X)\rightarrow \mathrm{Cl}(X)\rightarrow 
0.
\end{equation*}
where $T$ is the torus of the toric variety $X$.
The homomorphism $\alpha$ maps $m\in M$ to the $T$-invariant 
divisor $\mathrm{div}(\chi^{m})$, 
and the homomorphism $\mathrm{Div}_T(X)\rightarrow \mathrm{Cl}(X)$ maps every 
$T$-invariant Weil divisor on $X$ to its divisor class.
Here $M=\mathbb{Z}^3$ and
$\mathrm{Div}_T(X)\cong\mathbb{Z}^3$ since there are three rays in 
the fan of $X$.\\
By Lemma 3.2 and Lemma 4.3 above, we may assume the minimal 
generators of the cone $\sigma$ are of the form
$e_1,ae_1+be_2,ce_1+de_2+ee_3$, where 
$a,b,c,d,e\in\mathbb{Z}$ and $e\neq 0$.
Under the standard basis of $\mathbb{Z}^3$, the matrix 
representing $\alpha$ has the minimal generators of the cone 
$\sigma$ as its rows.
This matrix is
$A=\begin{pmatrix}
1 & 0 & 0\\
a & b & 0\\
c & d & e\\
\end{pmatrix}$.\\
Subtract $a$ times of the top row and add it to the second 
row, followed by subtracting $c$ times of the top row and 
adding it to the third row, we may assume $a=0=c$.
Hence, the Smith Normal form of the matrix $A$ is 
$\begin{pmatrix}
1 & 0 & 0\\
0 & b & 0\\
0 & 0 & e\\
\end{pmatrix}$.
Since the Smith Normal form is obtained from elementary row 
and column operations over the integers, which preserve 
determinants up to a sign, we have $|be|=|\delta|$, where
$\delta$ is the determinant of the matrix taking the minimal 
generators of the cone $\sigma$ as its columns.
Hence, the cokernel  of $\alpha$ is $\mathrm{Cl}(X)=A^1(X)$.
It is is a finite abelian group of order $|\delta|$, which
is isomorphic to 
$\mathbb{Z}/b\mathbb{Z}\oplus\mathbb{Z}/e\mathbb{Z}$.
\end{proof}
\section{Conjecture}
Here is a table obtained using the software SageMath.
It shows eight examples of 3-dimensional affine, simplicial 
toric varieties $X$ together with their Chow groups $A^1(X)$ 
and $A^2(X)$:\par
\begin{tabular}{ |p{3.9cm}|p{1cm}|p{2cm}|p{2cm}| }
\hline
minimal generators& $\delta$ & $A^1(X)$ & $A^2(X)$\\
\hline 
(1,0,0),(1,2,0),(1,2,4) &8& $C2\times C4$ & 0\\
\hline
(1,0,0),(1,3,0),(1,3,9) &27& $C3\times C9$ & 0\\
\hline
(1,0,0),(2,3,0),(3,5,7)& 21 & $C21$ & $C7$\\
\hline
(1,0,0),(2,5,0),(3,7,9)& 45& $C45$ & $C9$\\
\hline
(1,0,0),(5,7,11),(7,8,19)& 45& $C45$ & $C5$\\
\hline
(1,0,0),(3,5,0),(7,9,13)& 65& $C65$& $C13$\\
\hline
(1,0,0),(3,7,0),(5,8,11)& 77& $C77$&0\\
\hline
(1,0,0),(5,7,0),(7,8,19)& 133& $C133$& $C19$\\
\hline
\end{tabular}
\par
Let $X$ be any affine, simplicial toric variety over a field 
$k$.
Based on Theorem 4.4 and the computations above, we 
conjecture that the Chow group $A^1(X)$ has order 
$|\delta|$, where 
$\delta$ is the determinant of the matrix 
taking the minimal generators of the unique maximal cone of 
the fan of $X$ as its columns.
And the order of the Chow group $A^2(X)$ divides 
$|\delta|$.\par

\end{document}